# Asymptotics of the Eigenvalues of Two-Diagonal Jacobi Matrices

R. V. Kozhan

## 1. INTRODUCTION

In the Hilbert space $H := \ell_2(\mathbb{N})$ endowed with an inner product denoted by $(\,\cdot\,|\,\cdot\,)$, we consider the Jacobi operator

$$J = SA + AS^*, \tag{1.1}$$

where $S$ is the right-shift operator and $A$ is the diagonal operator with diagonal $(a_n)_{n=1}^\infty$, i.e.,

$$Se_n = e_{n+1}, \quad Ae_n = a_n e_n, \quad n \in \mathbb{N}.$$

Here $(e_n)_{n=1}^\infty$ is the standard basis in $\ell_2(\mathbb{N})$. It is assumed that the numerical sequence $(a_n)_{n=1}^\infty$ consists of positive numbers for which the following condition is satisfied:

$$\lim_{n \to \infty} \frac{a_{n+1}}{a_n} = 0. \tag{1.2}$$

The operator $J$ is compact and self-adjoint, and all of its eigenvalues are simple and nonzero. Moreover, its point spectrum is symmetric with respect to zero, i.e., $\sigma_p(J) = \{\pm \lambda_n \mid n \in \mathbb{N}\}$, where the sequence $(\lambda_n)_{n=1}^\infty$ consists of positive numbers that monotonically tend to zero.

The goal of this paper is to study the asymptotic behavior of the sequence $(\lambda_n)_{n=1}^\infty$. This problem was posed in connection with [1], where the asymptotic behavior of the sequence $(\lambda_n)_{n=1}^\infty$ was studied for the case in which $a_n = q^{n^s}$ (where $q \in (0,1)$, $s > 0$) and it was proved that the elements $\lambda_n$ satisfy the condition $a_{2n-1} > \lambda_n > a_{2n+1}$. The main result of this paper is as follows.

**Theorem 1.1.** *If condition* (1.2) *is satisfied, then*

$$\lambda_n = (1 + o(1))a_{2n-1}, \quad n \to \infty. \tag{1.3}$$

Note that our approach is based on the study of the quadratic form of the operator $J$ and is more elementary than the approach in [1], which used the apparatus of the theory of entire functions.

## 2. PROOF OF THEOREM 1.1

Suppose that $P_n \colon H \to H$ is an orthogonal projection operator on the subspace

$$H_n = \{x \in H \mid x = (x_k)_{k=1}^{\infty} \ \forall k > n : \ x_k = 0\}, \qquad n \in \mathbb{N}.$$

Let

$$J_n := P_n J P_n.$$

Also, let $T_n$ be the restriction of the operator $J_n$ to the space $H_n$.

**Lemma 2.1.** *The spectra of the operators $J$ and $J_n$, $n \in \mathbb{N}$, are symmetric with respect to zero. All the eigenvalues of the operator $J$ and all the nonzero eigenvalues of the operators $J_n$ are simple. Besides, the number $\lambda = 0$ is not an eigenvalue of the operator $J$.*

**Proof.** We denote by $U$ the diagonal operator that acts by the formula

$$U e_n = (-1)^n e_n, \qquad n \in \mathbb{N}.$$

Since $U^{-1} S U = -S$ and $U^{-1} A U = A$, we have

$$U^{-1} J U = -J, \qquad U^{-1} J_n U = -J_n;$$

hence it follows that the spectra of the operators $J$ and $J_n$ are symmetric. The proof of the other assertions of the lemma is obvious. □

We denote by $(\lambda_k)_{k=1}^{\infty}$ (respectively, by $(\lambda_k^{(n)})_{k=1}^{\infty}$) the sequence of nonnegative eigenvalues of the operator $J$ (respectively, of $J_n$) numbered in decreasing order, counting multiplicity. By $(s_k(X))_{k=1}^{\infty}$ we denote the sequence of $s$-numbers of a compact operator $X \colon H \to H$ (see [2]).

**Lemma 2.2.** *For arbitrary $n \in \mathbb{N}$, the following relation holds:*

$$\prod_{k=1}^{n} \lambda_k^{(2n)} = \prod_{k=1}^{n} a_{2k-1}.$$

**Proof.** Taking Lemma 2.1 into account, we find

$$\prod_{k=1}^{n} \lambda_k^{(2n)} = |\det T_{2n}|^{1/2}.$$

For arbitrary $n > 1$, simple calculations yield

$$\det T_{2n} = -a_{2n-1}^2 \det T_{2n-2},$$

and thus the lemma is proved. □

**Lemma 2.3.** *For sufficiently large $n_0 \in \mathbb{N}$, the following inequalities hold:*

$$\|J_m - J_n\| \leq a_n + a_{n+1}, \qquad m > n \geq n_0, \tag{2.1}$$

$$\|J - J_n\| \leq a_n + a_{n+1}, \qquad n \geq n_0. \tag{2.2}$$

**Proof.** It follows from condition (1.2) that $a_{n+1} \leq a_n$, beginning with some $n_0 \in \mathbb{N}$ and, therefore, for $n > n_0$ and $\|x\| = 1$, we obtain

$$\left|((J - J_n)x|x)\right| \leq 2\left|\sum_{k=n}^{\infty} a_k x_k \overline{x_{k+1}}\right| \leq \sum_{k=n}^{\infty} a_k\left(|x_k|^2 + |x_{k+1}|^2\right)$$

$$= \sum_{i=0}^{\infty} a_{n+2i}\left(|x_{n+2i}|^2 + |x_{n+2i+1}|^2\right) + \sum_{i=0}^{\infty} a_{n+2i+1}\left(|x_{n+2i+1}|^2 + |x_{n+2i+2}|^2\right)$$

$$\leq a_n + a_{n+1}. \tag{2.3}$$

Thus, inequality (2.2) is proved. Similarly, we can prove estimate (2.1). □

**Lemma 2.4.** *The inequality*

$$\lambda_k^{(2k)} \geq \frac{a_{2k-1}}{2}$$

*holds for infinitely many numbers $k \in \mathbb{N}$.*

**Proof.** Let us prove this by contradiction: suppose that there exists an $n_1 \in \mathbb{N}$ such that

$$\lambda_k^{(2k)} < \frac{a_{2k-1}}{2} \qquad \forall k \geq n_1.$$

By condition (1.2), there exists a number $n_2 \in \mathbb{N}$ such that $a_{n+1} \leq a_n/8$ for all $n \geq n_2$. Take $p = \max\{n_1, n_2\}$, and let $N > p$. Taking Lemma 2.2 into account, we obtain

$$\prod_{i=1}^{N} a_{2i-1} = \prod_{i=1}^{N} \lambda_i^{(2N)} = \prod_{i=1}^{p} \lambda_i^{(2N)} \cdot \prod_{i=p+1}^{N} \lambda_i^{(2N)}. \tag{2.4}$$

Since

$$|\lambda_i^{(2N)} - \lambda_i^{(2i)}| = |s_{2i}(J_{2n}) - s_{2i}(J_{2i})| \leq \|J_{2N} - J_{2i}\|,$$

in view of Lemma 2.3, we have

$$\lambda_i^{(2N)} < \lambda_i^{(2i)} + a_{2i} + a_{2i+1} \leq \frac{a_{2i-1}}{2} + 2 \cdot \frac{a_{2i-1}}{8} = \frac{3a_{2i-1}}{4}, \qquad p < i \leq N.$$

Since

$$|\lambda_i^{(2N)}| \leq \|J_{2N}\| = \|P_{2N} J P_{2N}\| \leq \|J\|, \qquad 1 \leq i \leq p,$$

from (2.4) we obtain

$$\prod_{i=1}^{N} a_{2i-1} \leq \|J\|^p \prod_{i=p+1}^{N} \frac{3a_{2i-1}}{4}$$

and, therefore,

$$\prod_{i=1}^{p} a_{2i-1} \leq \|J\|^p \left(\frac{3}{4}\right)^{N-p}.$$

Now, letting $N \to \infty$, we find that there are zeros among the numbers $a_n$. The resulting contradiction brings the proof to an end. □

**Proof of Theorem 1.1.** Choose an arbitrary $\varepsilon \in (0, 1/8)$. Then there exists (see (1.2)) a number $n_1 \in \mathbb{N}$ such that $a_n \leq \varepsilon a_{n-1}$ for all $n \geq n_1$.

First, let us show that for large $n \in \mathbb{N}$ the following inequality holds:

$$\lambda_n \leq (1+\varepsilon)a_{2n-1}. \tag{2.5}$$

Indeed, by inequality (2.2), there exists a number $n_0 > n_1$ for which

$$\|J - J_{2n-1}\| \leq a_{2n-1} + a_{2n}, \qquad n > n_0.$$

We can easily see that $\lambda_n = s_{2n}(J)$ and $s_{2n}(J_{2n-1}) = 0$. Therefore, for $n > n_0$ we have

$$\lambda_n = s_{2n}(J) - s_{2n}(J_{2n-1}) \leq \|J - J_{2n-1}\| \leq a_{2n-1} + a_{2n} \leq (1+\varepsilon)a_{2n-1}.$$

Next, we verify that, for sufficiently large $n \in \mathbb{N}$, the following inequality holds:

$$\lambda_n \geq (1 - 2\varepsilon)a_{2n-1}. \tag{2.6}$$

Using the result of Lemma 2.4, we choose a number $m > n_1$ so that $\lambda_m^{(2m)} \geq a_{2m-1}/2$. Since

$$\lambda_i^{(2m)} \geq \lambda_m^{(2m)} \geq \frac{a_{2m-1}}{2}, \qquad i = 1, \ldots, m,$$

by the spectral theorem, we have

$$(J_{2m}x | J_{2m}x) \geq \frac{1}{4}a_{2m-1}^2 \sum_{k=1}^{2m} |x_k|^2, \qquad x \in H_{2m}. \tag{2.7}$$

For an arbitrary $r > m+1$, consider the space

$$G_{2r}^m := \{(x_n)_{n=1}^\infty \in H_{2r} \mid x_{2m+1} = 0\}.$$

Note that $\dim G_{2r}^m = 2r - 1$. For arbitrary $x \in G_{2r}^m$, by direct calculations we obtain

$$(Jx|Jx) = (J_{2m}x|J_{2m}x) + a_{2m}^2|x_{2m}|^2$$
$$+ \sum_{k=2m+2}^{2r}(a_k^2 + a_{k-1}^2)|x_k|^2 + 2\operatorname{Re}\left(\sum_{k=2m}^{2r-2} a_k a_{k+1} x_k \overline{x_{k+2}}\right). \tag{2.8}$$

Since $a_k < \varepsilon a_{k-1}$ for all $k \geq 2m$, we have

$$2|a_k a_{k+1} x_k \overline{x_{k+2}}| \leq 2\varepsilon |a_{k-1} x_k| \cdot |a_{k+1} x_{k+2}| \leq \varepsilon(a_{k-1}^2|x_k|^2 + a_{k+1}^2|x_{k+2}|^2), \qquad k \geq 2m,$$

and, therefore,

$$2\operatorname{Re}\left(\sum_{k=2m}^{2r-2} a_k a_{k+1} x_k \overline{x_{k+2}}\right) \geq -\varepsilon a_{2m-1}^2 |x_{2m}|^2 - 2\varepsilon \sum_{k=2m+2}^{2r} a_{k-1}^2 |x_k|^2.$$

Hence, from (2.8) and (2.7), we find

$$(Jx|Jx) \geq \frac{a_{2m-1}^2}{4}\sum_{k=1}^{2m}|x_k|^2 - \varepsilon a_{2m-1}^2 |x_{2m}|^2 + (1-2\varepsilon)\sum_{k=2m+2}^{2r} a_{k-1}^2 |x_k|^2. \tag{2.9}$$

Since $\varepsilon \in (0, 1/8)$ and $r > m+1$, we have

$$a_{2m-1}^2\left(\frac{1}{4} - \varepsilon\right) \geq a_{2r-1}^2;$$

hence, taking (2.9) into account, we can write

$$(Jx|Jx) \geq (1-2\varepsilon)a_{2r-1}^2 \sum_{k=1}^{2r} |x_k|^2, \qquad x \in G_{2r}^m. \tag{2.10}$$

Now, applying the minimax principle [3, Chap. 8], in view of Lemma 2.1 and the relation $\dim G_{2r}^m = 2r - 1$, we obtain

$$\lambda_r \geq (1-2\varepsilon)^{1/2} a_{2r-1} \geq (1-2\varepsilon) a_{2r-1}, \qquad r > m+1. \tag{2.11}$$

Thus, inequality (2.6) is proved. Now, Theorem 1.1 readily follows from (2.5) and (2.6). □

I. FRANKO LVOV NATIONAL UNIVERSITY
*E-mail*: rkozhan@yahoo.com